\documentclass[11pt,oneside]{article}
\usepackage{amssymb,amsmath,amsthm}
\usepackage{url}
\usepackage[english]{babel}
\usepackage{pstricks,pst-node}
\renewcommand{\le}{\leqslant}
\renewcommand{\ge}{\geqslant}

\newcommand{\diag}{\mathrm{diag}}

\newcommand{\ZZ}{\mathbb{Z}}

\newcommand{\NN}{\mathbb{N}}

\newcommand{\LLL}{\mathcal{L}}

\newcommand{\SSS}{\mathcal{S}}

\newcommand{\vx}{\mathbf{x}}

\newcommand{\vy}{\mathbf{y}}

\newcommand{\va}{\mathbf{a}}

\newtheorem{theorem}{Theorem}
\newtheorem{proposition}{Proposition}
\newtheorem{problem}{Problem}

\newtheorem*{theoremB}{Theorem Bad}

\newcommand{\vA}{\mathbf{A}}

\newcommand{\contt}{K^{(t)}}
\newcommand{\contmn}{K^{(-1)}}

\begin{document}

\title{Continuant Diophantine equations}

\author{
 Dzmitry Badziahin
% \footnote{Research partially supported by EPSRC  Grant EP/E061613/1}
}

\maketitle
\begin{abstract}
We investigate a family of Diophantine polynomial equations which
involve continuant functions. In particular, given a polynomial
$P(x)\in \ZZ[x]$ and $n\in\NN$,  we consider the equation
$P(K_n(x_1,\ldots, x_n)) =
K_{n+1}(x_0,\ldots,x_n)K_{n+1}(x_1,\ldots, x_{n+1})$. We show that
with certain restrictions on $P(x)$ the set of its solutions has a
rich structure. In particular, we provide several ways of generating
new solutions from the existing ones.

In the last section we discuss the relation between the solutions of
the above Diophantine equation for arbitrary values of $n$ and
factorisations $P(m) = d_1d_2$ for integers $m,d_1$ and $d_2$.
\end{abstract}

\section{Introduction}

We start by introducing the generalised continuants $K_n^{(t)}$.
Given $t\in\ZZ\backslash \{0\}$, let $K_n^{(t)}$ be a polynomial of
$n$ variables which is defined by the recurrent formula
$$
K_0^{(t)}():=1;\quad K_1^{(t)}(x_1):=x_1;
$$
\begin{equation}\label{cont_eq}
\contt_{n+1}(x_1,\ldots,x_{n+1}) =
x_{n+1}\contt_n(x_1,\ldots,x_n)+t\contt_{n-1}(x_1,\ldots,x_{n-1}).
\end{equation}
To make the notation shorter, we denote
$$
\vx_{m,n}:= (x_m,x_{m+1},\ldots,x_n),\quad m\le n+1
$$
and
$$
\overline{\vx}_{m,n}:= (x_n,x_{n-1},\ldots,x_m),\quad m\le n+1.
$$
For $t=1$ the polynomial $\contt_n(\vx_{1,n})$ is the standard
continuant $K_n(\vx_{1,n})$ which plays an important role in the
theory of continued fractions. Indeed, a finite continued fraction
$[a_0;a_1,\ldots, a_n]$ is a rational number which enumerator and
denominator are $K_{n+1}(\va_{0,n})$ and $K_n(\va_{1,n})$
respectively. We refer the reader to~\cite[Chapter
X]{Hardy_Wright_2008} for details.

Fix $n\in \NN$, $n\ge 2$ and $t\in\ZZ\backslash \{0\}$. Let $P(x) =
c_0+c_1x+\ldots +c_dx^d$ be a polynomial with integer coefficients
such that
\begin{equation}\label{cond_p}
c_0=(-t)^n\quad\mbox{and}\quad x^d P((-t)^{n-1}x^{-1}) = C \cdot
P(x),
\end{equation}
where $C$ is some non-zero constant. One can easily check that monic
symmetric and antisymmetric polynomials are covered by this
property. Indeed, for those polynomials the condition~\eqref{cond_p}
is satisfied with $t=-1$. Many such polynomials also
satisfy~\eqref{cond_p} with $t=1$ and either an even or an odd value
of $n$.

The central Diophantine equation of this paper is
\begin{equation}\label{main_eq}
P(K^{(t)}_{n-1}(\vx_{1,n-1})) =
K^{(t)}_n(\vx_{0,n-1})K^{(t)}_n(\vx_{1,n}).
\end{equation}
Since $\deg(P)$ can be arbitrarily high, so can be the degree of
this equation. Despite this, we will show that it usually has
infinitely many solutions and moreover the set of solutions has
quite an interesting structure. In Section~\ref{sec3} we show that
for $t=\pm 1$ every solution $\vx_{0,n}$ of~\eqref{main_eq}
generates a sequence $\SSS = \ldots,x_{-1},x_0,\ldots,x_n,\ldots$ of
integers such that any tuple $\vx_{m,m+n}$ is also a solution
of~\eqref{main_eq}. In Section~\ref{sec4} we provide a different way
of constructing infinitely many solutions $\vx$ of the
equation~\eqref{main_eq}. We do it by noticing that
$\contt_n(\vx_{0,n-1})=\pm 1$ implies that there exists $x_n\in\ZZ$
such that $\vx_{0,n}$ is a solution of~\eqref{main_eq}. We believe
that in many cases different solutions $\vx_{0,n}$ achieved in this
way generate different sequences $\SSS$. We do not prove this
statement formally, however we can see this in examples considered
in the paper. Next, in Section~\ref{sec5} we show that, base on a
solution of~\eqref{main_eq} for a particular value of $n$, we can
construct new solutions of the equation for the same polynomial
$P(x)$ and bigger values of $n$. That gives us the third method of
generating solutions of~\eqref{main_eq}. In the last section we
notice the relation between the factorisations $P(m) = d_1d_2$ for
integer values of $m,d_1, d_2$ and the solutions of~\eqref{main_eq}
where $t=1$, $P(x)$ is fixed and $n$ may vary. At the end, by
considering factorisations of the polynomial $m^4+1$ for small
integer values $m$ we conclude that in general the three methods for
generating new solutions of~\eqref{main_eq} described above are
still not sufficient to provide all solutions of this equation.
Therefore the problem of classifying all the solutions
of~\eqref{main_eq} still waits for its discoverer.

To the best of authors knowledge, equations of the
form~\eqref{main_eq} are mostly uncovered in the literature. We can
only refer to a paper~\cite{Badziahin_2016} where the case $n=2$ and
some particular types of the polynomials $P(x)$ were considered. On
the other hand, continuants play an important role in the solutions
of some Diophantine equations. The most classical example is Pell's
equation $x^2-dy^2=\pm 1$ where $d$ is not a perfect square. It is
well known that all its solutions $(x,y)$ are
$(K_{kn}(\va_{0,kn-1}), K_{kn-1}(\va_{1,kn-1}))$, where
$[a_0;a_1,a_2,\ldots]$ is a continued fraction of $\sqrt{d}$, $n$ is
the length of its period and $k$ is any positive integer. We refer
to \cite[Section IV.11]{Davenport_1982} for details. Also Schinzel
recently used continuants to find all solutions of the equation
$x^2+x+1 = yz$~\cite{Schinzel_2015}.

\section{Properties of generalised continuants}

Continuants naturally arise from the following matrix identity which
can easily be checked from the recurrent formula~\eqref{cont_eq}:
$$
\left(\begin{array}{cc} 0&t\\
1&x_1
\end{array}\right)\cdot\left(\begin{array}{cc} 0&t\\
1&x_2
\end{array}\right)\cdots\left(\begin{array}{cc} 0&t\\
1&x_n
\end{array}\right) = \left(\begin{array}{cc} t\contt_{n-2}(\vx_{2,n-1})&t\contt_{n-1}(\vx_{2,n})\\
\contt_{n-1}(\vx_{1,n-1})& \contt_n(\vx_{1,n})
\end{array}\right).
$$
It in turn provides the following properties of continuants which we
will use later.
\begin{equation}\label{cont_prop1}
\contt_n(\vx_{1,n})\contt_{n-2}(\vx_{2,n-1}) -
\contt_{n-1}(\vx_{1,n-1})\contt_{n-1}(\vx_{2,n}) = (-1)^nt^{n-1}.
\end{equation}
\begin{equation}\label{cont_prop2}
\begin{array}{rl}\contt_{n+1}(\vx_{1,n+1}) &= x_n\contt_n(\vx_{1,n}) +
t\contt_{n-1}(\vx_{1,n-1})\\ &= x_1 \contt_n(\vx_{2,n+1}) +
t\contt_{n-1}(\vx_{3,n+1}).
\end{array}
\end{equation}
\begin{equation}\label{cont_prop3}
\contt_n(\vx_{1,n}) = \contt_n(\overline{\vx}_{1,n}).
\end{equation}
The first two properties are trivial:~\eqref{cont_prop1} can be
derived by considering the determinant of both hand sides
and~\eqref{cont_prop2} is an application of the standard induction.
To check~\eqref{cont_prop3} one notices that
$$
\left(\begin{array}{cc} 0&t\\
1&x_1
\end{array}\right)^T = \left(\begin{array}{cc} t&0\\
0&1
\end{array}\right)^{-1}\cdot\left(\begin{array}{cc} 0&t\\
1&x_1
\end{array}\right)\cdot \left(\begin{array}{cc} t&0\\
0&1
\end{array}\right).
$$
Then, by transposing the both sides of the matrix equation,
multiplying by $\diag(t,1)$ from the left and by $\diag(t,1)^{-1}$
from the right we get
$$
\left(\begin{array}{cc} 0&t\\
1&x_n
\end{array}\right)\cdot\left(\begin{array}{cc} 0&t\\
1&x_{n-1}
\end{array}\right)\cdots\left(\begin{array}{cc} 0&t\\
1&x_1
\end{array}\right) = \left(\begin{array}{cc} t\contt_{n-2}(\vx_{2,n-1})&t\contt_{n-1}(\vx_{1,n-1})\\
\contt_{n-1}(\vx_{2,n})& \contt_n(\vx_{1,n})
\end{array}\right).
$$
Now property~\eqref{cont_prop3} follows straightforwardly.

In the end of this section we provide another property of
generalised continuants. From~\eqref{cont_prop2} we can derive that
\begin{equation}\label{cont_prop4}
\contt_n(\vx_{1,n}) = \contt_{n+1}(1,x_1-t,\vx_{2,n}) =
\contt_{n+2}(\vx_{1,n-1},x_n-t,1).
\end{equation}

\section{Sequences of equation
\eqref{main_eq} solutions}\label{sec3}

Property~\eqref{cont_prop3} already reveals some symmetry of the
solutions of~\eqref{main_eq}: if $\vx_{0,n}$ solves it then so does
its ``inverse'' $\overline{\vx}_{0,n}$. However one can say much
more about their structure. The next theorem, which is one of the
main results of this paper, says that every solution
of~\eqref{main_eq} generates a chain of solutions.

\begin{theorem}\label{th1}
Let $P(x)$ be a polynomial with integer coefficients, which
satisfies~\eqref{cond_p}. Assume that $\vx_{0,n} \in \ZZ^{n+1}$ is
an integer solution of the equation~\eqref{main_eq} such that each
term $x_0,\ldots, x_n$ is coprime with $t$. Then there exists
$x_{n+1}\in \ZZ$ and $x_{-1}\in \ZZ$ such that $\vx_{1,n+1}$ and
$\vx_{-1,n-1}$ are also solutions of~\eqref{main_eq}.

Moreover the value $x_{n+1}$ (respectively $x_{-1}$) is uniquely
defined by $\vx_{0,n}$ unless $\contt_{n-1}(\vx_{2,n})=0$
(respectively $\contt_{n-1}(\vx_{0,n-2})=0$). In the latter case,
$\vx_{1,n+1}$ is a solution of~\eqref{main_eq} for any integer
$x_{n+1}$.
\end{theorem}
\proof Firstly note that coprimeness conditions and recurrent
formula~\eqref{cont_eq} imply that for any $m,k\in\ZZ_{\ge 0}$ such
that $0\le m\le n-k+1$ the values $\contt_{k}(\vx_{m,m+k-1})$ are
coprime with $t$. Moreover,~\eqref{cont_prop1} implies that
$\contt_{n-1}(\vx_{2,n})$ and $\contt_n(\vx_{1,n})$ are coprime.
Assume that $\contt_{n-1}(\vx_{2,n})\neq 0$. Then we can
apply~\eqref{cont_prop1} to get
$$
\contt_{n-1}(\vx_{1,n-1}) \equiv (-t)^{n-1}
(\contt_{n-1}(\vx_{2,n}))^{-1}\pmod {\contt_n(\vx_{1,n})}.
$$
Substituting this into the equation~\eqref{main_eq} gives us
$$
0\equiv P(\contt_{n-1}(\vx_{1,n-1})) \equiv P((-t)^{n-1}
(\contt_{n-1}(\vx_{2,n}))^{-1}) \stackrel{\eqref{cond_p}}{\equiv}
P(\contt_{n-1}(\vx_{2,n}))\pmod {\contt_n(\vx_{1,n})}.
$$
Therefore there is an integer $T$ such that
\begin{equation}\label{eq_th1}
P(\contt_{n-1}(\vx_{2,n})) = T\cdot \contt_n(\vx_{1,n}).
\end{equation}
Look at equation~\eqref{eq_th1} modulo $\contt_{n-1}(\vx_{2,n})$
(recall that $\contt_{n-1}(\vx_{2,n})\neq 0$). By the first property
in~\eqref{cond_p} the left hand side of~\eqref{eq_th1} is congruent
to $(-t)^n$. On the other hand,~\eqref{cont_prop2} implies that its
right hand side is congruent to $tT\cdot \contt_{n-2}(\vx_{3,n})$.
Hence we have the following congruence
$$
T\cdot \contt_{n-2}(\vx_{3,n}) \equiv (-t)^{n-2}\cdot t \pmod{
\contt_{n-1}(\vx_{2,n})}.
$$
Aditionally,~\eqref{cont_prop1} implies that
$$
\contt_{n-2}(\vx_{2,n-1})\contt_{n-2}(\vx_{3,n}) \equiv (-t)^{n-2}
\pmod{\contt_{n-1}(\vx_{2,n})}.
$$
By comparing the last two congruences and by cancelling
$\contt_{n-2}(\vx_{3,n})$ we get $T\equiv
t\contt_{n-2}(\vx_{2,n-1})$ $\pmod{\contt_{n-1}(\vx_{2,n})}$ and
therefore there exists $x_{n+1}\in \ZZ$ such that
$$
T = x_{n+1}\contt_{n-1}(\vx_{2,n})+
t\contt_{n-2}(\vx_{2,n-1})\quad\mbox{or}
$$
\begin{equation}\label{eq_t}
T = \contt_n(\vx_{2,n+1}).
\end{equation}
Note that, since $\contt_{n-1}(\vx_{2,n})\neq 0$, $x_{n+1}$ is
defined uniquely (as a solution of linear equation).

Now suppose the contrary, $\contt_{n-1}(\vx_{2,n})= 0$. Then the
value $\contt_n(\vx_{1,n})$ equals either 1 or $-1$ because
otherwise it is not coprime with $\contt_{n-1}(\vx_{2,n})=0$. This
fact implies~\eqref{eq_th1} for some $T\in\ZZ$ because $\pm1$
divides every integer number. By the first property
in~\eqref{cond_p} the left hand side of~\eqref{eq_th1} equals
$(-t)^n$. On the other hand by~\eqref{cont_prop2} the right hand
side equals $tT\cdot \contt_{n-2}(\vx_{3,n})$. Property
\eqref{cont_prop1} implies that
$\contt_{n-2}(\vx_{2,n-1})\contt_{n-2}(\vx_{3,n})=(-t)^{n-2}$.
Finally we get $T = t\contt_{n-2}(\vx_{2,n-1})$. Since
$\contt_{n-1}(\vx_{2,n})=0$, the equation~\eqref{eq_t} is satisfied
for any integer $x_{n+1}$.

To conclude, we have shown that for any integer solution $\vx_{0,n}$
 of~\eqref{main_eq} there exists  $x_{n+1}\in
\ZZ$ such that $\vx_{1,n+1}$ is also a solution of~\eqref{main_eq}.
Analogous arguments allow us to find $x_{-1}\in\ZZ$ such that
$\vx_{-1,n-1}$ is a solution of~\eqref{main_eq}.
\endproof

Note that from Theorem~\ref{th1} we can extract the formula for
$x_{n+1}$:
\begin{equation}\label{eq_an}
x_{n+1} = \frac{P(\contt_{n-1}(\vx_{2,n})) -
t\contt_{n-2}(\vx_{2,n-1})\contt_n(\vx_{1,n})}{\contt_{n-1}(\vx_{2,n})\contt_n(\vx_{1,n})}.
\end{equation}
If $\contt_{n-1}(\vx_{2,n})=0$ then any $x_{n+1}$ gives a solution
of~\eqref{main_eq}. Finally if $\contt_{n}(\vx_{1,n})=0$ then we
necessarily have $P(\contt_{n-1}(\vx_{2,n}))=0$ and
$$
x_{n+1} =
-\frac{t\contt_{n-2}(\vx_{2,n-1})}{\contt_{n-1}(\vx_{2,n})}.
$$

%Theorem~\ref{th1} tells that any solution $\vx_{0,n}$ of the
%equation~\eqref{main_eq} comes with the sequence $\SSS$
%$$
%\ldots,x_{-2},x_{-1},x_0,\ldots, x_n,x_{n+1},\ldots
%$$
%of integer numbers such that any $n+1$ consecutive elements of this
%sequence form a solution of~\eqref{main_eq}. It can be finite or
%infinite in any of two directions. It terminates at the position $m$
%on the right side as soon as it encounters
%$\contt_{n-1}(\vx_{m-n+2,m})=0$. The same is for the left side of
%the sequence: it terminates at the position $m$ as soon as
%$\contt_{n-1}(\vx_{m,m+n-2})=0$. We denote the set of all
%$(n+1)$-tuples of consecutive numbers from this sequence by
%$\LLL(\vx_{0,n})$. So each set $\LLL(\vx_{0,n})$ gives us the series
%of solutions of~\eqref{main_eq}. As we will see later, in many cases
%these sets are infinite.

We can slightly modify Theorem~\ref{th1} to show that the solution
$\vx_{0,n}$ is usually uniquely determined by its $n$ first
elements. In fact, we can prove

%
%set $\LLL(\vx_{0,n})$ can be uniquely determined by $n$ consecutive
%elements in $\SSS$ instead of $n+1$ provided that these $n$ elements
%are not at the endpoints of the sequence $\SSS$.

\begin{theorem}\label{th2}
Let $\vx_{0,n-1}$ be an integer $n$-tuple such that all its terms
are coprime with $t$. If the value $\contt_n(\vx_{0,n-1})$ divides
$P(\contt_{n-1}(\vx_{1,n-1}))$ then there exists $x_n\in\ZZ$ such
that $\vx_{0,n}$ is a solution of~\eqref{main_eq}. Moreover, this
value is unique, provided that $\contt_{n-1}(\vx_{1,n-1})\neq 0$.
\end{theorem}

\proof We have that there exists an integer $T$ such that
$$
P(\contt_{n-1}(\vx_{1,n-1})) = T\cdot \contt_n(\vx_{0,n-1}).
$$
Then we repeat the arguments of Theorem~\ref{th1} starting
from~\eqref{eq_th1} with $\vx_{0,n-1}$ in place of $\vx_{1,n}$.
\endproof

As before we can compute the value of $x_n$ from Theorem~\ref{th2}.
In this case it can be computed by formula~\eqref{eq_an} with all
the scripts of $x$'s shifted by one unit left.

Note that in both Theorems we assumed that all elements of
$\vx_{0,n}$ (respectively of $\vx_{0,n-1}$ in Theorem~\ref{th2}) are
coprime with $t$. Unfortunately, we can not guarantee the same
property for $x_{-1}$ and $x_{n+1}$ (respectively for $x_n$).
However in two cases, $t=\pm 1$, the coprimality conditions for
$\vx_{0,n}$ become trivial and can be removed. For convenience, in
further discussion we will always set~$t$ to one of these two
values.

Theorem~\ref{th1} shows that any solution $\vx_{0,n}$ of the
equation~\eqref{main_eq} comes with the sequence $\SSS$
$$
\ldots,x_{-2},x_{-1},x_0,\ldots, x_n,x_{n+1},\ldots
$$
of integer numbers such that any its $n+1$ consecutive elements form
a solution of~\eqref{main_eq}. It can be finite or infinite from
either side. It terminates at the position $m$ from the right if it
encounters $\contt_{n-1}(\vx_{m-n+2,m})=0$. The same is for the left
side of the sequence: it terminates at the position $m$ if
$\contt_{n-1}(\vx_{m,m+n-2})=0$. Furthermore, if $\vx_{0,n-1}$ is
not the utmost $n$-tuple of~$\SSS$ then this sequence is uniquely
defined by $\vx_{0,n-1}$. We denote the set of all $(n+1)$-tuples of
consecutive numbers from $\SSS$ by $\LLL(\vx_{0,n})$ and call it
{\em a chain}. As we will see later, in many cases chains
$\LLL(\vx_{0,n})$ are infinite and provide an infinite set of
different solutions of~\eqref{main_eq}. If
$\contt_{n-1}(\vx_{1,n-1})\neq 0$ we will also use the notation
$\LLL(\vx_{0,n-1})$ for the chain~$\LLL(\vx_{0,n})$ to emphasize
that it is uniquely defined by $\vx_{0,n-1}$.

As we mentioned before, Property~\eqref{cont_prop3} suggests that
the set of solutions of~\eqref{main_eq} is closed under inverting
the terms in $(n+1)$-tuples $\vx_{0,n}$. Moreover, careful
investigation of the formula~\eqref{eq_an} shows that the chain
$\LLL(\overline{\vx}_{0,n})$ consists of all inverted $(n+1)$-tuples
from $\LLL(\vx_{0,n})$. We denote such a chain by
$\overline{\LLL}(\vx_{0,n})$, i.e.
$$
\overline{\LLL}(\vx_{0,n}):= \LLL(\overline{\vx}_{0,n}).
$$

%Theorem~\ref{th1} shows that each solution $\vx_{0,n}$
%of~\eqref{main_eq} defines a sequence
%$\ldots,x_{-2},x_{-1},x_0,x_1\ldots$ such that every $n+1$
%consecutive elements of it are also solutions of~\eqref{main_eq}.
%The sequence can end up with $n-1$ elements $\vx_{k,k+n-2}$ such
%that $\contt_{n-1}(\vx_{k,k+n-2})=0$ from either side or can be
%infinite. Theorem~\ref{th2} states that $x_n$ is uniquely defined
%given the previous $n$ terms $\vx_{0,n-1}$. Also every $n$ elements
%of this sequence inverted also form a solution of~\eqref{main_eq}.
%We denote the set of solutions from this sequence together with
%inverted solutions by $\LLL(\vx_{0,n-1})$. More precisely,
%$$
%\LLL(\vx_{0,n-1}):= \{\vx_{k,k+n},\overline{\vx}_{k,k+n} \;:\;
%k_*\le k\le k^*\}
%$$
%where the bound $k_*$ is defined by the equation:
%$$
%\contt_{n-1}(\vx_{k_*,k_*+n-2}) = 0\;\mbox{ or }\;k_* =
%-\infty\;\mbox{ if }\;\contt_{n-1}(\vx_{k,k+n-2})\neq 0\;\mbox{ for
%any }\; k\le 0.
%$$
%The bound $k^*$ is defined analogously:
%$$
%\contt_{n-1}(\vx_{k^*+1,k^*+n-1}) = 0\;\mbox{ or }\;k^* =
%\infty\;\mbox{ if }\;\contt_{n-1}(\vx_{k+1,k+n-1})\neq 0\;\mbox{ for
%any }\;k\ge 0.
%$$
It is quite straightforward to check that if
$\vx\in\LLL(\vx_{0,n-1})$ then $\LLL(\vx) = \LLL(\vx_{0,n-1})$. In
other words the sequences $\LLL(\vx_{0,n-1})$ define the equivalence
classes for the roots of~\eqref{main_eq}: $\vx$ is equivalent to
$\vy$ if $\LLL(\vx) = \LLL(\vy)$.

We finish this section by the natural question:

\begin{problem}
Given $n\in\NN$, $t=\pm 1$ and a polynomial $P(x)$ which
satisfies~\eqref{cond_p}, classify all chains $\LLL(\vx_{0,n-1})$ of
solutions of~\eqref{main_eq}.
\end{problem}

In general this problem seems to be very difficult. The case $n=2$
and $t=1$ is considered in~\cite{Badziahin_2016}, where the complete
classification of chains $\LLL(x_0,x_1)$ such that
$x_0,x_1\in\ZZ_{\ge 0}$ for polynomials $P$ of degree four is
provided. For certain polynomials of higher
degrees~\cite{Badziahin_2016} gives only a partial classification of
chains. We will present the result for the polynomial $P(x) = x^4+1$
from that paper at the end of the next section.

\section{Equation $\contt_n(x_1,\ldots,x_n) = 1$}\label{sec4}

Recall that now we only consider the cases $t=\pm 1$. Additionally,
from now on we assume that the polynomial $P(x)$ is fixed and so is
the equation~\eqref{main_eq}. Therefore, whenever we mention the
chain $\LLL(\vx)$, we mean the chain of solutions of that fixed
equation.

Note that if $\contt_n(\vx_{0,n-1}) = 1$ then the conditions of
Theorem~\ref{th2} are definitely satisfied. This fact allows us to
construct many chains $\LLL(\vx_{0,n-1})$. Indeed, it is not
difficult to check that the equation $\contt_n(\vx_{0,n-1})=1$ has
infinitely many solutions for all $n\ge 2$ if $t=1$ and for all
$n\ge 3$ if $t=-1$. For example, the value $K^{(1)}_2(x_1,x_2) =
x_1x_2+1$ equals one if and only if one of the terms $x_1$ or $x_2$
equals zero. Therefore we have an infinite collection of chains
$\LLL(0,x_2)$ where $x_2$ is arbitrary. In general, any solution of
$\contt_n(x_1,\ldots, x_n) = 1$ gives a chain $\LLL(\vx_{1,n})$.
Therefore the classification of all solutions $\contt_n(\vx_{1,n}) =
1$ gives a family of chains $\LLL(\vx_{1,n})$. Unfortunately, as we
will see in Section~\ref{sec_fact} in the general case this family
does not give us a full classification of equation~\eqref{main_eq}
solutions.

%For example, Theorem~\ref{th2} provides us with solution
%of~\eqref{main_eq} as soon as $\contt_n(\vx_{0,n-1}) = 1$. Then for
%every such $n$-tuple Theorem~\ref{th1} gives us a series of
%solutions $\LLL(\vx_{0,n-1})$. Finally,

Lets firstly concentrate on the case $t=1$. A complete family of
solutions of Diophantine equations $K_n(x_1,\ldots,x_n)=1$ for small
values $n$ can be found with relatively standard elementary methods.
A quick search of the literature did not reveal any research on the
equations of this type. Therefore here we briefly consider these
equations for $n=2,3$ and 4 to explain general ideas behind their
solutions.

Because of the Property~\eqref{cont_prop3} we can restrict our
search of the solutions to the case $(x_1,\ldots,x_n)\le_l
(x_n,x_{n-1},\ldots,x_1)$ where by $\le_l$ we mean being smaller or
equal in the lexicographical order.

The case $n=2$ is straightforward:
\begin{proposition}
The solutions of
$$
K_2(x_1,x_2) = 1
$$
consist of one-parametric series $(0,a)$, $a\in\ZZ$ together with
its inverse $(a,0)$.
\end{proposition}

For higher values of $n$ the number of different families of
solutions grows quickly. However we still have a finite amount of
them.
%infinite series of solutions together with several extra
%solutions.

\begin{proposition}\label{prop3}
The equation
$$
K_3(x_1,x_2,x_3) = 1
$$
has
\begin{itemize}
\item three one-parameter series of solutions
$(0,a,1),(a,0,1-a),(1,-1,a)$, $a\in\ZZ$ together with their inverses
and
\item three extra solutions $(-3,1,-2), (-1,3,-1), (-2,2,-1)$ together
with their inverses.
\end{itemize}
\end{proposition}
\proof The equation can be rewritten as follows
$$
x_1x_2x_3 + x_1+x_3=1.
$$
Notice that at least one of the three values $|x_1|,|x_2|$ and
$|x_3|$ must be less than two. Indeed, otherwise we have that
$|x_1|,|x_3|$ are at most a quarter of $|x_1x_2x_3|$ and $1\le
1/8|x_1x_2x_3|$. Therefore
$$
|x_1x_2x_3| > |x_1|+|x_3|+1
$$
which contradicts the equation.

Then we consider several cases.
\begin{itemize}
\item $x_1=0$. This gives a series of the solutions $(0,a,1)$, $a\in\ZZ$.
\item $x_2=0$ gives solutions $(a,0,1-a)$, $a\in\ZZ$.
\item $x_3=0$ gives solutions $(1,a,0)$, $a\in\ZZ$.
\item $x_1=1$ gives the equation $x_3(x_2+1) = 0$. The case $x_3=0$
has already been considered therefore $x_2=-1$ which gives solutions
$(1,-1,a)$, $a\in\ZZ$.
\item $x_2=1$ gives the equation $x_1x_3+x_1+x_3=1$. After
rewriting it as $(x_1+1)(x_3+1) = 2$ and by taking all possible
factorisations of 2 we get new solutions $(-2,1,-3)$ and
$(-3,1,-2)$.
\item $x_3=1$ gives solutions $(a,-1,1)$, $a\in\ZZ$.
\item $x_1=-1$ gives the equation $x_3(1-x_2) = 2$. By taking all
possible factorisations of 2 we get extra solutions $(-1,3,-1)$ and
$(-1,2,-2)$.
\item $x_2=-1$ gives the equation $(x_1-1)(x_3-1)=0$ with no new
solutions.
\item $x_3=-1$ gives one more additional solution $(-2,2,-1)$.
\end{itemize}\endproof

In the equation $K_4(\vx_{1,4})=1$ for $n=4$ we can also prove that
at least one of the variables $|x_1|,\ldots,|x_4|$ is less than two.
Then it again suffices to consider the finite amount of cases. Here
we only provide the classification of the solutions of
$K_4(\vx_{1,4})=1$ and leave the rigorous proof to the interested
reader.

\begin{theorem}
In total, the equation $K_4(x_1,x_2,x_3,x_4)=1$ has
\begin{itemize}
\item three two-parameter series of solutions $(0,a,b,0), (0,a,0,b),
(a,0,-a,b)$, $a,b\in\ZZ$ together with their inverses;
\item three one-parameter series of solutions $(1,a,-1,1), (-1,a,1,-1),
(a,-1,1,a)$ together with their inverses;
\item 15 isolated solutions $(-4,1,-2,2), (-3,1,-3,1),
(-3,1,-2,3), (-3,2,-1,3), (-2,1,-4,1),$ $(-2,1,-3,2), (-2,2,-2,1),
(-2,2,-1,4), (-2,3,-1,2), (-1,-1,-1,1), (-1,2,-3,1),$ $(-1,2,-2,2),
(-1,3,-2,1), (-1,3,-1,3),(-1,4,-1,2)$ and their inverses.
\end{itemize}
\end{theorem}

As we can see, the amount of different series of roots of the
equation $K_n(x_1,\ldots,x_n)=1$ grows rapidly with $n$. We expect
that it will continue growing like that for bigger $n$ too.

For higher values of $n$ we only show that there are infinitely many
solutions of the equation~\eqref{eq_kn}:
\begin{proposition}\label{prop2}
Equation
\begin{equation}\label{eq_kn}
K_n(x_1,\ldots,x_n)=1
\end{equation}
has infinitely many integer
solutions for any $n\ge 2$.
\end{proposition}

\proof We have already showed this for $n\le 4$. Now assume that
$n\ge 5$. By using Property~\eqref{cont_prop2} of continuants we get
that $(x_1,\ldots,x_{n-1},0)$ is a solution of~\eqref{eq_kn} as soon
as $(x_1,\ldots,x_{n-2})$ is a solution of
$K_{n-2}(x_1,\ldots,x_{n-1})=1$. As we already know, the latter has
infinitely many solutions.
\endproof

Theorem~\ref{th2} and Proposition~\ref{prop2} together suggest the
way of finding infinitely many different
sequences~$\LLL(\vx_{0,n-1})$. We consider all $n$-tuples
$\vx_{0,n-1}$ such that $K_n(\vx_{0,n-1})=1$
(Proposition~\ref{prop2} says that there are infinitely many of
them) and then by Theorem~\ref{th2} we construct $x_n\in\ZZ$ such
that $\vx_{0,n}$ is a solution of~\eqref{main_eq}. Theoretically all
 $(n+1)$-tuples constructed in this way may belong to the same sequence
$\LLL(\vx_{0,n-1})$ but on practise this is usually not the case.

The case $t=-1$ can be considered in a similar way.

\begin{theorem}\label{th3}
The equation $\contmn_2(x_1,x_2) = 1$ has two solutions $(1,2)$ and
$(-2,-1)$ together with their inverses. The equation
$$
\contmn_3(x_1,x_2,x_3) = 1
$$
has
\begin{itemize}
\item three infinite series of solutions $(0,a,-1), (a,0,-1-a),
(-1,-1,a)$ together with their inverses.
\item four additional solutions $(1,2,2), (1,3,1), (2,1,3)$ together with their
inverses.
\end{itemize}
Finally for $n\ge 3$ the equation $\contmn_n(\vx_{1,n}) = 1$ has
infinitely many solutions.
\end{theorem}

\proof For $n=2$ the equation can be rewritten as $x_1x_2=2$. By
considering all factorisations of 2 we get the solutions $(-2,-1),
(-1,-2), (1,2), (2,1)$.

For $n=3$ we notice that $K_3^{(-1)}(x_1,x_2,x_3) = 1$ if and only
if $K_3(-x_1,x_2,-x_3) = 1$. All solutions of the latter equation
are provided by Proposition~\ref{prop3}.

Finally, by~\eqref{cont_prop2} we have that
$\contmn_{n+2}(\vx_{1,n},a,0) = -\contmn_n(\vx_{1,n})$ for any
$a\in\ZZ$. Therefore if the equation $\contmn_n(\vx_{1,n})=-1$ has
at least one solution then the equation $\contmn_{n+2}(\vx_{1,n+2})
= 1$ has infinitely many solutions. By analogous reason, if
$\contmn_{n-2}(\vx_{1,n-2})=1$ has at least one solution then
$\contmn_{n+2}(\vx_{1,n+2})$ has infinitely many solutions.

By the definition of the generalised continuants, the equation
$\contmn_n(\vx_{1,n})=1$ has a solution for $n=0$ and $n=1$,  we
already know that it has solutions for $n=2$ and $n=3$ and basic
induction finishes the proof for $n\ge 4$.
\endproof

To conclude, the solutions of $\contt_n(\vx_{1,n}) = 1$ generate the
chains $\LLL(\vx_{1,n})$. By the same reason the solutions of the
similar equation $\contt_n(\vx_{1,n}) = -1$ may generate more chains
$\LLL(\vx_{1,n})$. By this method we can construct at least partial
classification of chains. We call the chain $\LLL(\vx_{1,n})$
nonstandard if $\forall \vx\in \LLL(\vx_{1,n})$ one has
$\contt_n(\vx)\neq \pm 1$. Then the Problem~A can be speified a bit
further.

\begin{problem}
Given $n\in\NN$, $t=\pm1$ and a polynomial $P(x)$ which
satisfies~\eqref{cond_p}, classify all nonstandard chains
$\LLL(\vx_{1,n})$.
\end{problem}

\subsection{The case $n=2$ and $P(x)=x^4+1$}

Let $n=2$. The paper~\cite{Badziahin_2016} contains a full
classification of chains $\LLL(x_0,x_1)$ with nonnegative $x_0$ and
$x_1$ for polynomials $P(x)$ of degree 4. With some efforts that
classification can be extended to all chains $\LLL(x_0,x_1)$. We
demonstrate this statement for a model example $P(x)=x^4+1$.

In this case the result from~\cite{Badziahin_2016} states

\begin{theoremB}
The positive integer solutions of the equation
\begin{equation}\label{eq_thb}
x_1^4+1 = (x_0x_1+1)(x_1x_2+1)
\end{equation}
are elements of chains $\LLL(0,a)$ where $a\in\NN$.
\end{theoremB}

Firstly note that the solutions of~\eqref{eq_thb} are close under
changing the sign: $(x_0,x_1,x_2)$ is a solution of~\eqref{eq_thb}
if and only if $(-x_0,-x_1,-x_2)$ is. This fact immediately gives us
the classification of all negative solutions of~\eqref{eq_thb}: they
are elements of the chains $\LLL(0,a)$ where $a<0$.

Next, the solutions with at least one zero term can be easily found.
If $x_0=0$ then $\vx_{0,2}$ is an element of a chain $\LLL(0,a)$ for
some $a\in\ZZ$. Similarly solutions with $x_2=0$ are elements of
chains $\overline{\LLL}(0,a)$, $a\in\ZZ$. Finally, a straightforward
check shows that any triple $(x_0,0,x_2)$ is a solution
of~\eqref{eq_thb}.

The remaining case is when the signs of the solutions
$(x_0,x_1,x_2)$ alternate. Without loss of generality assume that
$x_0>0>x_1$. Then from~\eqref{eq_thb} we have that $x_2$ must also
be positive. The chain $\LLL(x_0,x_1)$ can not contain zero since we
already classified all chain with zeroes and none of them have terms
with alternating signs. Therefore, by consecutive inspection of
elements $x_3,x_4,\ldots$ from the chain $\LLL(x_0,x_1)$ we get that
terms of any solution $\vx_{k,k+2}\in\LLL(x_0,x_2)$ have alternating
signs. By replacing $\vx_{0,2}$ with its inverse
$\overline{\vx}_{0,2}$, if needed, we may guarantee that
$|x_0|>|x_1|$. Therefore the chain $\LLL(x_0,x_1)$ must contain an
element $\vx_{k,k+2}$ such that $|x_k|\ge|x_{k+1}|\le |x_{k+2}|$,
otherwise it would contain zero, which is impossible. By
rewriting~\eqref{eq_thb} with
$$
|x_{k+1}|^4+1 = (|x_kx_{k+1}|-1)(|x_kx_{k+1}|-1)
$$
we can easily classify all such integer solutions. Up to a sign
there is only one such a solution: $|x_k|=|x_{k+2}|=2$ and
$|x_{k+1}|=1$. Hence $\vx_{0,2}$ is an element of either
$\LLL(2,-1,2)$ or $\LLL(-2,1,-2)$.

As a conclusion we end this section with the complete classification
of integer solutions of~\eqref{eq_thb}.

\begin{theorem}\label{th5}
Every integer solution $(x_0,x_1,x_2)$ of~\eqref{eq_thb} falls into
one of the following categories:
\begin{itemize}
\item $\vx_{0,2}$ is an element of $\LLL(0,a)$ or
$\overline{\LLL}(0,a)$ where $a\in\ZZ\backslash\{0\}$;
\item $\vx_{0,2}$ has $x_1=0$;
\item $\vx_{0,2}$ is an element of $\LLL(-2,1,-2)$ or
$\LLL(2,-1,2)$.
\end{itemize}
\end{theorem}

A trivial inspection of all three cases above shows that for $n=2$
ad $P(x) = x^4+1$ all chains $\LLL(x_0,x_1)$ are standard.

\section{Other ways of generating solutions
of~\eqref{main_eq}}\label{sec5}

In this section we will show that each solution $\vx_{0,n+1}$ of the
equation~\eqref{main_eq} generates many other solutions
of~\eqref{main_eq} for the same polynomial $P(x)$ but for different
values of $n$. Fix a polynomial $P(x)$ and consider
Condition~\eqref{cond_p} for various values $n$. Notice that for
$t=1$, if $P(x)$ satisfies~\eqref{cond_p} for some $n$ then it
satisfies the same condition for all values of $n$ of the same
parity. Furthermore, for $t=-1$ if $P(x)$ satisfies~\eqref{cond_p}
for some $n$ then it satisfies the same condition for all $n$. Based
on this observation we denote by $\vA^t(P)$ the set of all integer
solutions of~\eqref{main_eq} for a given $t$ and $P(x)$ and for any
integer $n$ such that $P(x)$ satisfies~\eqref{cond_p}.

Further in this section we will always assume that $\vx_{0,n+1}$ is
a solutions of~\eqref{main_eq}. To start with, assume that $t=1$. We
have $K_{n+1}(\vx_{0,n})\mid P(K_n(\vx_{1,n}))$. Then we also have
that for any $d\in \ZZ$,
$$
K_{n+3}(0,a,\vx_{0,n}) = K_{n+1}(\vx_{0,n}) \mid
P(K_n(\vx_{1,n})+aK_{n+1}(\vx_{0,n})) = P(K_{n+2}(a,\vx_{0,n})).
$$
By applying Theorem~\ref{th2} we find the value $x'_{n+1}$ such that
$(0,a,\vx_{0,n},x'_{n+1})$ is a new solution of~\eqref{main_eq} for
the value $n+2$. Since the parity of $n$ and $n+2$ coincide, we have
that $(0,a,\vx_{0,n},x'_{n+1})\in \vA^1(P)$. Moreover if
$K_n(\vx_{1,n})\neq 0$ then the value $x'_{n+1}$ is unique and can
be computed by formula~\eqref{eq_an}. Otherwise
$(0,a,\vx_{0,n},x'_{n+1})$ is a solution of the main equation for
any integer $x'_{n+1}$. Therefore we can define the maps from
$\vA^1(P)$ to itself:
$$
f_a\;:\;\{ \vx_{0,n+1}\in \vA^1(P)\;:\; K_n(\vx_{1,n})\neq 0\} \to
\vA^1(P)
$$
such that $f_a(\vx_{0,n+1}) = (0,a,\vx_{0,n},x'_{n+1})$. Also define
$$
f^*_{a,b}\;:\;\{\vx_{0,n+1}\in \vA^1(P)\;:\; K_n(\vx_{1,n})= 0\} \to
\vA^1(P)
$$
by $f_a(\vx_{0,n+1}) = (0,a,\vx_{0,n},b)$. These two functions $f_a$
and $f^*_{a,b}$ provide a way of generating new solutions
to~\eqref{main_eq} in addition to chains $\LLL(\vx_{0,n})$.

From now on Let's assume that $P(x)$ is either even or odd (i.e.
$P(-x) = \pm P(x)$). Then the same generating procedure will work
for $t=-1$ too. Indeed,
$$
\contmn_{n+3}(0,a,\vx_{0,n}) = -\contmn_{n+1}(\vx_{0,n}) \mid
P(-\contmn_n(\vx_{1,n})+a\contmn_{n+1}(\vx_{0,n})) =
P(\contmn_{n+2}(a,\vx_{0,n})).
$$
Therefore in that case we can define two maps $f^{(-1)}_a$ and
$f^{*(-1)}_{a,b}$ for $t=-1$ in the same way as $f_a$ and
$f^*_{a,b}$.

Finally, we provide one more way of generating new elements of
$\vA^t(P)$ for even and odd polynomials $P$. By~\eqref{cont_prop4}
we have that $\contt_{n+2}(1,x_0-t,\vx_{1,n}) =
\contt_{n+1}(\vx_{0,n})$. On the other hand, by
Condition~\eqref{cont_prop2} we have
$$
\contt_{n+1}(x_0-t,\vx_{1,n}) = \contt_{n+1}(\vx_{0,n}) -
t\contt_n(\vx_{1,n})
$$
Therefore for $t=\pm 1$ a divisibility $\contt_{n+1}(\vx_{0,n})\mid
P(\contt_n(\vx_{1,n}))$ implies
$$
\contt_{n+2}(1,x_0-t,\vx_{1,n}) \mid
P(\contt_{n+1}(x_0-t,\vx_{1,n})).
$$
For $t=-1$ this gives a new solution of the
equation~\eqref{main_eq}. Define a new map
$$
g\;:\; \{\vx_{0,n}\in \ZZ^{n+1}\;:\; n\in \ZZ_{\ge 0}\}\to
\{\vx_{0,n}\in \ZZ^{n+1}\;:\; n\in \ZZ_{\ge 0}\}
$$
by $g(\vx_{0,n}) = (1,x_0-t,\vx_{1,n})$. Notice that in the case
$x_0=1$ the inverse map $g^{-1}$ is correctly defined. For $t=-1$,
$g$ maps $\vA^{-1}(P)$ to itself. However, since $g$ changes the
parity of the length of $\vx_{0,n}$, $g(\vA^1(P))$ is not
necessarily in $\vA^1(P)$. On the other hand it can be fixed by
introducing one more map
$$
h\;:\; \{\vx_{0,n}\in \ZZ^{n+1}\;:\; n\in \ZZ_{\ge 0}\}\to
\{\vx_{0,n}\in \ZZ^{n+1}\;:\; n\in \ZZ_{\ge 0}\}
$$
by $h(\vx_{0,n}) = (\vx_{0,n-1},x_n-t,1)$. Notice that
$$
\contt_{n+2}(\vx_{0,n-1},x_n-t,1) = \contt_{n+1}(\vx_{0,n})\mid
P(\contt_n(\vx_{1,n})) = P(\contt_{n+1}(\vx_{1,n},x_{n+1}-t)).
$$
Again, as before, function $h$ maps $\vA^{-1}(P)$ to itself. On the
other hand it also changes the parity of the size of vectors $\vx$.
Finally notice that $g\circ h = h\circ g$ does not change the parity
of $n$, therefore $g\circ h$ maps $\vA^1(P)$ to itself. Finally, if
correctly defined, the maps $g^{\pm 1}\circ h^{\pm 1}$ also map
$\vA^1(P)$ to itself and therefore provide another way to generate
new solutions of~\eqref{main_eq}.

We conclude this section with the table which contains all maps
discussed here. Given $\vx_{0,n}\in \vA^t(P)$ we can produce the
following new elements of $\vx_{0,n}\in \vA^t(P)$:

$$
\begin{array}{|c|c|}
\hline t=1&t=-1\\
\hline f_a(\vx_{0,n})\quad\mbox{if}\; K_n(\vx_{1,n})\neq 0&
f^{(-1)}_a(\vx_{0,n})\quad\mbox{if}\; K^{-1}_n(\vx_{1,n})\neq 0\\
\hline f_{a,b}(\vx_{0,n})\quad\mbox{if}\; K_n(\vx_{1,n})= 0&
f^{(-1)}_{a,b}(\vx_{0,n})\quad\mbox{if}\; K^{-1}_n(\vx_{1,n})=0\\
\hline g\circ h(\vx_{0,n}) & g(\vx_{0,n})\\
\hline g\circ h^{-1} (\vx_{0,n})\quad \mbox{if }\;x_n=1&
h(\vx_{0,n})\\
\hline g^{-1}\circ h (\vx_{0,n})\quad \mbox{if }\;x_0=1&
g^{-1}(\vx_{0,n})\quad\mbox{if }\; x_0=1\\
\hline g^{-1}\circ h^{-1} (\vx_{0,n})\quad \mbox{if }\;x_0=x_n=1&
h^{-1}(\vx_{0,n})\quad\mbox{if }\; x_n=1\\
\hline
\end{array}
$$

\section{Relation with the factorisation of values
$P(m)$}\label{sec_fact}

Fix $t=1$ and a polynomial $P(x)$ which satisfies~\eqref{cond_p} for
some integer value $n_0$. In the previous section we observed that
$P(x)$ then satisfies conditions~\eqref{cond_p} for all $n$ of the
same parity as $n_0$.

Consider an arbitrary factorisation $P(m) = d_1d_2$ of the value
$P(m)$ at some positive integer value $m$ into the product of two
integer factors. Write $d_1/m$ as a continued fraction:
$$
\frac{d_1}{m} = [a_0;a_1,\ldots, a_{n-1}]
$$
where $a_0\in \ZZ$ and $a_1,\ldots, a_n\in\NN$. Moreover, we can
always choose such repesentation that $n - n_0$ is even. Since
$t=1$, conditions~\eqref{cond_p} imply that $m$ and $P(m)$ are
coprime, therefore $d_1/m$ is irreducible continued fraction. From
the theory of continued fractions we have
\begin{equation}\label{eq_md1}
m = K_{n-1}(\va_{1,n-1})\quad\mbox{and}\quad d_1 =K_n(\va_{0,n-1})
\end{equation}
Therefore the $n$-tuple $\va_{0,n-1}$ satisfies the conditions of
Theorem~\ref{th2}. Moreover, $K_{n-1}(\va_{1,n-1}) = m\neq 0$ which
implies that there exists unique value $a_n$ such that $\va_{0,n}$
is a solution of~\eqref{main_eq}.

To conclude, any factorisation of $P(m)$ for positive integer $m$
generates a solution $\va_{0,n}$ of~\eqref{main_eq} with
\begin{equation}\label{cond_va}
a_0,a_n\in\ZZ\quad\mbox{and}\quad a_1,\ldots, a_{n-1}\in\NN.
\end{equation}
Conversely, for any solution $\va_{0,n}$ which
satisfies~\eqref{cond_va} we have
$$
P(K_{n-1}(\va_{1,n-1})) = K_n(\va_{0,n-1})\cdot K_n(\va_{1,n}).
$$
Denote $m = K_{n-1}(\va_{1,n-1})\in\NN$, $d_1 = K_n(\va_{0,n-1})$
and $d_2 = K_n(\va_{1,n})$. Whence we have a bijection between all
factorisations of $P(m)$ and the solutions $\va_{0,n}$
of~\eqref{main_eq} satisfying~\eqref{cond_va}. Therefore if we
understand the full structure of the set $\vA^1(P)$ then we can
generate factorisations of $P(m)$ in a constructive way. In
particular, we will be able to construct values $P(m)$ which
factorisation as a product of primes satisfies the prescribed
properties. This is quite important for cryptosystems like RSA.

\subsection{Example}

For better understanding of the equation~\eqref{main_eq} and its
link to the factorisations of $P(m)$ we look at our model example
$P(x) = x^4+1$ and $t=1$. It is easy to check that $P(x)$
satisfies~\eqref{cond_p} for any even value $n$.

As discussed in the previous chapters we have a number of ways to
create new solutions in $\vA^1(P)$ based on already known ones. For
any given solution $\vx_{0,n}$ we can create the chain
$\LLL(\vx_{0,n})$ of solutions, we can use maps $f_a$ and
$f^*_{a,b}$, finally we have maps $g^{\pm 1}\circ h^{\pm 1}$. The
natural question is then as follows: {\it given all these maps and
finitely many solutions of~\eqref{main_eq}, can we generate all the
set $\vA^1(P)$?} We do not know the formal answer to this question,
however, as we will see in a moment, it seems to be negative. In any
case, with help of the described maps we can construct a subclass of
$\vA^1(P)$ and therefore the subclass of the factorisations of
numbers like $m^4+1$.

We start with the table which covers the factorisations of $m^4+1$
for small values of $m$. In the ``$\vx$ and a sequence $\SSS$''
column we write the solution $\vx$ in brackets and then we write the
sequence $\SSS$ around it.

$$
\begin{array}{|c|c|c|c|}
\hline
m&\mbox{factorisation}&\vx \mbox{ and a sequence $\SSS$}&\mbox{notation}\\
\hline 1&1^4+1 = 1\cdot 2& {\bf (0,1,1)},0&\vx^{(1)}\\
\hline 2&2^4+1 = 1\cdot 17& {\bf (0,2,8)},30,112,\ldots& \vx^{(2)}\\
\hline 3&3^4+1 = 1\cdot 82& {\bf (0,3,27)},240,2133,\ldots& \vx^{(3)}\\
\hline 3&3^4+1 = 2\cdot 41& -1,{\bf (0,1,1,1,13)},480,23422307,\ldots& \vx^{(4)} = f_1(\vx^{(1)})\\
\hline 4&4^4+1 = 1\cdot 257& {\bf (0,4,64)},1020,\ldots& \vx^{(5)}\\
\hline 5&5^4+1 = 1\cdot 626& {\bf (0,5,125)},3120,\ldots& \vx^{(6)}\\
\hline 5&5^4+1 = 2\cdot 313& -2,{\bf (0,2,1,1,62)},6240,\ldots& \vx^{(7)} = f_2(\vx^{(1)})\\
\hline 6&6^4+1 = 1\cdot 1297& {\bf (0,6,216)},7770,\ldots& \vx^{(8)}\\
\hline 7&7^4+1 = 1\cdot 2402& {\bf (0,7,343)},16800,\ldots& \vx^{(9)}\\
\hline 7&7^4+1 = 2\cdot 1201& -3,{\bf (0,3,1,1,171)},33600,\ldots& \vx^{(10)}= f_3(\vx^{(1)})\\
\hline 8&8^4+1 = 1\cdot 4097& {\bf (0,8,512)},32760,\ldots& \vx^{(11)}\\
\hline 8&8^4+1 = 17\cdot 241& 0,{\bf (2,8,30)},112,\ldots& \vx^{(12)} \in \LLL(\vx^{(2)})\\
\hline 9&9^4+1 = 1\cdot 6562& {\bf (0,9,729)},59040,\ldots& \vx^{(13)}\\
\hline 9&9^4+1 = 2\cdot 3281& -4,{\bf (0,4,1,1,364)},118080,\ldots& \vx^{(14)} = f_4(\vx^{(1)})\\
\hline 9&9^4+1 = 17\cdot 386& \ldots,21011,198,{\bf (1,1,7,1,42)},104543,\ldots& \vx^{(15)} = g\circ h(\vx^{(12)})\\
\hline 9&9^4+1 = 34\cdot 193& \ldots,42022,99,{\bf (3,1,3,2,21)},17487,\ldots& \vx^{(16)}, \mbox{new?}\\
\hline 10&10^4+1 = 1\cdot 10001& {\bf (0,10,1000)},99990,\ldots& \vx^{(17)}\\
\hline 10&10^4+1 = 73\cdot 137& \ldots, 1817, {\bf (7,3,2,1,13)},503,\ldots& \vx^{(18)}, \mbox{new?}\\
\hline
\end{array}
$$

Based on this table, we make several observations. As we can see,
all the factorisations of $m^4+1$ for $1\le m\le 10$, except two,
are generated by solutions $(0,a,a^3)\in\vA^1(P)$ which in turn come
from the solutions of $K_2(x_1,x_2) = 1$. However the remaining two
products, namely $9^4 + 1 = 34\cdot 193$ and $10^4+1 = 73\cdot 137$,
seem to generate new elements of $\vA^1(P)$ which do not intersect
with those generated by $(0,a,a^3)$. By considering larger values of
$m$ one can reveal more such elements. For example $\vx =
(7,2,2,2,19)\in \vA^1(P)$, which comes from the factorisation $12^4
+ 1 = 89\cdot 233$, is another example of this type. We believe that
infinitely many generators $\vx$ on top of $(0,a,a^3)$ are required
to generate the whole set $\vA^1(P)$. It will be interesting to see
the formal proof of this statement.

Let's consider $\vx^{(4)}$ and look at several first factorisations
which are linked with elements of the chain $\LLL(\vx^{(4)})$:
$$\begin{array}{lcl}
\vx_{0,4} = (0,1,1,1,13)&:& 3^4 = 2\cdot 41\\
\vx_{1,5} = (1,1,1,13,480)&:& 27^4 + 1 = 41\cdot 12962\\
\vx_{2,6} = (1,1,13,480,23422307)&:& 6721^4 + 1 = 12962\cdot
157421325361\\
\end{array}$$
The next solution $\vx_{3,7}$ gives a factorisation of the number
$(K_3(13,480,23422307))^4 + 1 = 146178618000^4 + 1$. This figures
show that elements of $\vx_{n,n+4}$ from $\LLL(\vx)$ provide
non-trivial factorisations of values $m^4+1$ and that the numbers
$m$ grow very quickly as $n$ grows. With help of the
formula~\eqref{eq_an} we can estimate the rate of growth of $x_n$.
Assuming that $0<x_1\le x_2\le x_3\le x_4$ we have that
$K_4(\vx_{1,4}) =
x_1K_3(\vx_{2,4})+K_2(\vx_{3,4})<(x_1+1)K_3(\vx_{2,4})$, therefore
$$
x_5 = \frac{(K_{3}(\vx_{2,4}))^4+1 -
K_{2}(\vx_{2,3})K_4(\vx_{1,4})}{K_{3}(\vx_{2,4})K_4(\vx_{1,4})}\ge
\frac{K_3(\vx_{2,4})^2}{x_1+1}-1.
$$
So definitely $x_5>\frac12 x_4^2x_3^2x_2$.

Theorem~\ref{th5} states that for $n=2$ and $P(x) = x^4+1$ all
chains $\LLL(x_0,x_1)$ are standard. However the table suggests that
for the same polynomial $P(x)$ and for higher values $n$ nonstandard
chains do exist. For example, consider $\vx^{(16)}=(3,1,3,2,21)$.
Indeed, one can show that for $n\ge 4$ elements $x_n$ in the chain
are strictly increasing. Also for $n<0$ elements $x_n$ strictly
decrease. We leave the rigorous proof of this statement to the
reader. Therefore a quick inspection shows that $K_4(\vx_{m,m+3})>
1$ for any four consecutive elements from the chain
$\LLL(3,1,3,2,21)$.

%
%\section{Applications}
%
%\subsection{$n=2$}
%
%Firstly consider the case $n=2$. Then Equation~\eqref{main_eq}
%transforms to
%\begin{equation}\label{main_n2}
%P(a_1) = (a_0a_1+1)(a_1a_2+1).
%\end{equation}
%The condition $x^dP(x^{-d}) = P(x)$ is equivalent to the polynomial
%$P$ being symmetric. This case was closely considered in the
%previous paper of the author~\cite{badziahin_2015} for quartic $P$
%and for $P(x) = x^8 + 1$. In particular it was shown there that
%\begin{itemize}
%\item For $P(x) = x^4+1$ all non-negative solutions solutions of~\eqref{main_n2}
%lie in the sequences $\LLL(0,n,n^3)$. Since $K_2(0,n) = 1$, all of
%these sequences are of the form provided at the end of the previous
%section.
%\item For an arbitrary symmetric quartic polynomial $P(x)$ the
%algorithm is provided for finding all sequences $\LLL(a_0,a_1,a_2)$
%of solutions. There may be several infinite collections of them
%together with the finite amount of exceptional sequences. At least
%not all of them necessarily contain a pair $a_n,a_{n+1}$ with
%$K_2(a_n,a_{n+1}) = 1$.
%\item Finally for $P(x) = x^8+1$ three infinite series of sequences
%were found:
%\begin{enumerate}
%\item $\LLL(0,n,n^7)$,\quad $n\in \NN$;
%\item $\LLL(n^5,n^3,n^5(n^8-1))$,\quad %n\in\NN$;
%\item $\LLL(\frac12t(d-1),\frac12t(d+1))$ where $t,d\in\ZZ_{\ge 0}$ satisfy
%the equation $d^2-2t^2=-7$.
%\end{enumerate}
%\end{itemize}

\end{document}